\pdfoutput=1
\RequirePackage{ifpdf}
\ifpdf 
\documentclass[pdftex]{sigma}
\else
\documentclass{sigma}
\fi

\numberwithin{equation}{section}

\newtheorem{Theorem}{Theorem}[section]
\newtheorem{Corollary}[Theorem]{Corollary}
\newtheorem{Lemma}[Theorem]{Lemma}
 { \theoremstyle{definition}
\newtheorem{Definition}[Theorem]{Definition}
\newtheorem{Example}[Theorem]{Example}
\newtheorem{Remark}[Theorem]{Remark} }

\usepackage{enumerate} 
\usepackage{xargs} 

\newcommand{\R}{\mathbb{R}} 
\renewcommand{\b}[1]{\mathbf{#1}} 
\newcommand{\set}[2]{\left\{ #1 \,\middle|\, #2 \right\}} 
\newcommand{\ssm}{\smallsetminus} 

\newcommand{\polygon}{\mathrm{P}} 
\newcommand{\dissection}{\mathrm{D}} 
\newcommand{\accordionComplex}{\mathcal{AC}} 
\newcommand{\sign}[3]{\varepsilon \big( {#1} \in {#2}\,|\,{#3} \big)} 
\newcommand{\SSS}{\reflectbox{$\mathsf{Z}$}} 
\newcommand{\ZZZ}{\mathsf{Z}} 
\newcommand{\VVV}{{\mathsf{{V \hspace{-.1686cm} I\,}}}} 
\newcommand{\gvector}[2]{\mathbf{g}(#1 \,|\, #2)} 
\newcommandx{\gvectorFan}[1][1=\dissection_\circ]{\mathcal{F}^\mathbf{g}(#1)} 
\newcommand{\complexP}{\mathcal{P}} 
\newcommand{\complexQ}{\mathcal{Q}} 
\newcommand{\quiver}{\mathrm{Q}} 
\newcommand{\relations}{\mathrm{I}} 
\newcommand{\stau}{\operatorname{s} \! \tau \! \operatorname{-tilt}}
\newcommand{\cpx}{2 \! \operatorname{-cpx}}
\newcommand{\silt}{2 \! \operatorname{-silt}}
\newcommand{\siltingComplex}{\mathcal{SC}}
\newcommand{\tiltingComplex}{\operatorname{s} \! \tau\mathcal{C}}
\newcommand{\MOD}{\mbox{{\rm mod \!}}}
\newcommand{\proj}{\operatorname{{\rm proj }}}
\newcommand{\Hom}[1]{\operatorname{{\rm Hom}}_{#1}}
\newcommand{\End}[1]{\operatorname{\rm End}_{#1}}

\begin{document}

\newcommand{\arXivNumber}{1710.02119}

\renewcommand{\PaperNumber}{045}

\FirstPageHeading

\ShortArticleName{A $\tau$-Tilting Approach to Dissections of Polygons}

\ArticleName{A $\boldsymbol{\tau}$-Tilting Approach to Dissections of Polygons}

\Author{Vincent PILAUD~$^a$, Pierre-Guy PLAMONDON~$^b$ and Salvatore STELLA~$^c$}

\AuthorNameForHeading{V.~Pilaud, P.-G.~Plamondon and S.~Stella}

\Address{$^a$~CNRS \& LIX, \'Ecole Polytechnique, Palaiseau, France}
\EmailD{\href{mailto:vincent.pilaud@lix.polytechnique.fr}{vincent.pilaud@lix.polytechnique.fr}} \URLaddressD{\url{http://www.lix.polytechnique.fr/~pilaud/}}

\Address{$^b$~Laboratoire de Math\'ematiques d'Orsay, Universit\'e Paris-Sud,\\ \hphantom{$^b$}~CNRS, Universit\'e Paris-Saclay, France}
\EmailD{\href{mailto:pierre-guy.plamondon@math.u-psud.fr}{pierre-guy.plamondon@math.u-psud.fr}} \URLaddressD{\url{https://www.math.u-psud.fr/~plamondon/}}

\Address{$^c$~University of Haifa, Israel}
\EmailD{\href{mailto:stella@math.haifa.ac.il}{stella@math.haifa.ac.il}}
\URLaddressD{\url{http://math.haifa.ac.il/~stella/}}

\ArticleDates{Received February 26, 2018, in final form May 10, 2018; Published online May 14, 2018}

\Abstract{We show that any accordion complex associated to a dissection of a convex polygon is isomorphic to the support $\tau$-tilting simplicial complex of an explicit finite dimensional algebra. To this end, we prove a property of some induced subcomplexes of support $\tau$-tilting simplicial complexes of finite dimensional algebras.}


\Classification{16G10; 16G20; 05E10}

\section{Introduction}

\begin{figure}[t]\centering
\includegraphics[width=.88\textwidth]{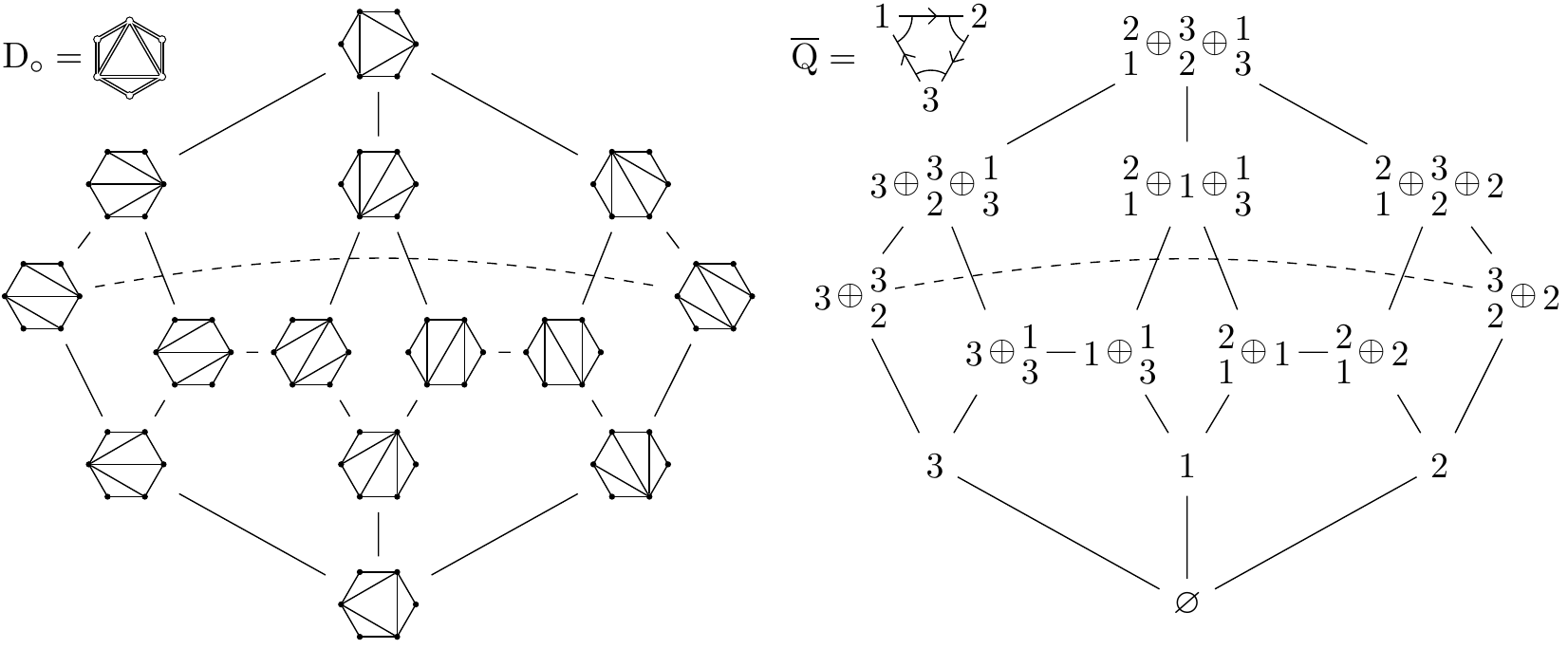}
\caption{The exchange graph on triangulations of a hexagon (left) and the exchange graph on support $\tau$-tilting modules of the quiver with relations~$\overline{\quiver}$ (right).}\label{fig:associahedron}
\end{figure}

The theory of cluster algebras gave rise to several interpretations of associahedra~\cite{ StasheffI,StasheffII, Tamari}. Fig.~\ref{fig:associahedron} shows two such interpretations for the rank~$3$ associahedron: as the exchange graph of triangulations of a hexagon and as the exchange graph of support $\tau$-tilting modules over the cluster tilted algebra whose quiver with relations is as depicted. This follows from results in the setting of the ``additive categorification of cluster algebras'' that was initiated in \cite{BuanMarshReinekeReitenTodorov, CalderoChapotonSchiffler}.

\begin{figure}[t]\centering
\includegraphics[width=.88\textwidth]{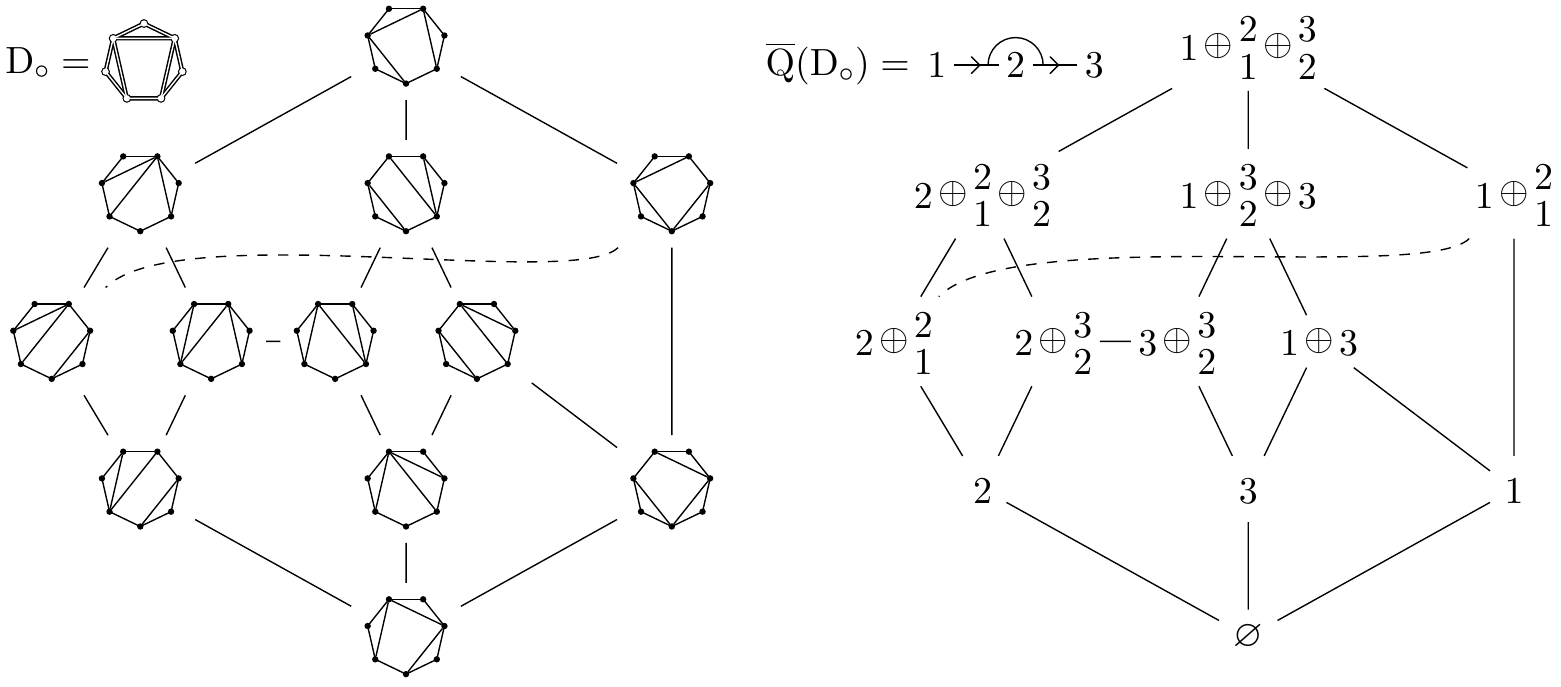}
\caption{The $\dissection_\circ$-accordion complex of the dissection~$\dissection_\circ$ of Fig.~\ref{fig:exmAccordionDissections} (left) and the $2$-term silting complex of the quiver~$\overline{\quiver}(\dissection_\circ)$ (right).}\label{fig:accordiohedron}
\end{figure}

F.~Chapoton observed a similar isomorphism between the exchange graph of certain dissections of a heptagon and that of support $\tau$-tilting modules over the path algebra of the quiver \includegraphics[scale=.35]{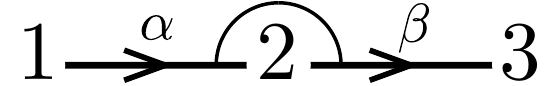}, subject to the relation $\beta\alpha = 0$. Fig.~\ref{fig:accordiohedron} shows these two exchange graphs, which can be found in~\cite[Fig.~7]{Chapoton-quadrangulations} and in~\cite[Example~6.4]{AdachiIyamaReiten}.

The purpose of this note is to show that this isomorphism is an avatar of a more general result in the theory of $\tau$-tilting modules. Any basic finite dimensional algebra $\Lambda$ gives rise to an exchange graph on support $\tau$-tilting modules. This exchange graph is the dual graph of the support $\tau$-tilting complex~\cite{AdachiIyamaReiten} (see Section~\ref{sec:recollections}). Let~$\{e_1, \dots, e_n\}$ be a complete set of primitive pairwise orthogonal idempotents of~$\Lambda$. Let~$J$ be a non-empty subset of~$[n]$ and~$e_J := \sum_{j \in J} e_j$. The following result forms the algebraic core of the paper (see Section~\ref{sec:recollections} for definitions and Section~\ref{sec:algRes} for the proof).

\begin{Theorem}\label{thm:mainAlgThm}The support $\tau$-tilting complex of~$e_J \Lambda e_J$ is isomorphic to the subcomplex of the support $\tau$-tilting complex of~$\Lambda$ induced by the support $\tau$-tilting modules whose $\b{g}$-vectors' coordinates vanish outside of~$J$.
\end{Theorem}

With this algebraic result, we can explain the isomorphism of Fig.~\ref{fig:accordiohedron}, and extend it to dissections of any polygon. Any reference dissection of a polygon gives rise to an exchange graph on certain dissections. This exchange graph is the dual graph of the accordion complex studied in~\cite{Chapoton-quadrangulations, GarverMcConville, MannevillePilaud-accordion} (see Section~\ref{sec:accordions}).

\begin{Theorem}\label{thm:mainDissection}Any accordion complex is isomorphic to the support $\tau$-tilting simplicial complex of an explicit finite dimensional algebra. Thus, the corresponding exchange graphs are isomorphic.
\end{Theorem}

Note that this statement could be seen as a consequence of a more general result~\cite[Proposition~2.44]{PaluPilaudPlamondon}. However, the latter combines several isomorphisms (from dissections, via non-kissing paths, through non-attracting chords, to support $\tau$-tilting modules) and the proof is quite intricate. In the present paper, we will easily deduce Theorem~\ref{thm:mainDissection} from Theorem~\ref{thm:mainAlgThm} and from the known case of triangulations of polygons (see Section~\ref{sec:accordions}).

The explicit finite dimensional algebras that appear in Theorem~\ref{thm:mainDissection} had previously appeared in~\cite{DavidRoesler, DavidRoeslerSchiffler}, where gentle algebras are associated to dissections of any surface without puncture. Applying Theorem~\ref{thm:mainAlgThm} to these algebras, one could obtain an analogue of the accordion complex for these surfaces.

\section[Recollections on $\tau$-tilting theory]{Recollections on $\boldsymbol{\tau}$-tilting theory}\label{sec:recollections}

The theory of $\tau$-tilting modules was introduced in~\cite{AdachiIyamaReiten}, and we mainly follow this source. Let $k$ be an algebraically closed field, let $\Lambda$ be a basic finite-dimensional $k$-algebra, and let $\{e_1, \dots, e_n\}$ be a complete set of pairwise orthogonal idempotents in~$\Lambda$. Denote by $\MOD \Lambda$ the category of finite-dimensional right $\Lambda$-modules, and by $\proj \Lambda$ its full subcategory of projective modules. We denote by $\tau$ the Auslander--Reiten translation of $\MOD \Lambda$ (see, for instance, \cite[Chapter~IV]{AssemSimsonSkowronski}). For any $\Lambda$-module $M$, we denote by $|M|$ the number of pairwise non-isomorphic direct summands appearing in any decomposition of $M$ into indecomposable modules.

\subsection[Support $\tau$-tilting pairs]{Support $\boldsymbol{\tau}$-tilting pairs}

Following~\cite[Definition~0.1]{AdachiIyamaReiten}, we say that a $\Lambda$-module $M$ is
\begin{itemize}\itemsep=0pt
 \item \textit{$\tau$-rigid} if $\Hom{\Lambda}(M, \tau M) = 0$;
 \item \textit{$\tau$-tilting} if it is $\tau$-rigid and $|M|=|\Lambda|$;
 \item \textit{support $\tau$-tilting} if there exists an idempotent $e$ of $\Lambda$ such that $e$ is in the annihilator of~$M$ and~$M$ is a~$\tau$-tilting $\Lambda/(e)$-module.
\end{itemize}
Support $\tau$-tilting modules always exist: $\Lambda$ itself and the zero module are two examples.

It is useful to keep track of the idempotents in the annihilator of a support $\tau$-tilting module. For this reason, we will follow~\cite[Definition~0.3]{AdachiIyamaReiten} and call a pair $(M,P)$, with $M\in \MOD \Lambda$ and $P\in \proj \Lambda$, a
\begin{itemize}\itemsep=0pt
 \item \textit{$\tau$-rigid pair} if $M$ is $\tau$-rigid and $\Hom{\Lambda}(P,M) = 0$;
 \item \textit{support $\tau$-tilting pair} if it is a $\tau$-rigid pair and $|M| + |P| = |\Lambda|$;
 \item \textit{almost complete support $\tau$-tilting pair} if it is a $\tau$-rigid pair and $|M| + |P| = |\Lambda|-1$.
\end{itemize}
We will say that the pair $(M,P)$ is \textit{basic} if both $M$ and $P$ are basic $\Lambda$-modules. We define direct sums of pairs componentwise.

One of the main theorems of~\cite{AdachiIyamaReiten} is the following.

\begin{Theorem}[{\cite[Theorem~0.4]{AdachiIyamaReiten}}] A basic almost complete support $\tau$-tilting pair is a direct summand of exactly two basic support $\tau$-tilting pairs.
\end{Theorem}

\begin{Definition}The \textit{support $\tau$-tilting complex} of~$\Lambda$ is the simplicial complex~$\tiltingComplex(\Lambda)$ whose vertices are the isomorphism classes of indecomposable $\tau$-rigid pairs and whose faces are sets of $\tau$-rigid pairs whose direct sum is rigid. The \textit{exchange graph}~$\stau(\Lambda)$ is the dual graph of~$\tiltingComplex(\Lambda)$, i.e., the graph whose vertices are isomorphism classes of basic support $\tau$-tilting pairs, and where two vertices are joined by an edge whenever the corresponding support $\tau$-tilting pairs differ by exactly one direct summand.
\end{Definition}

\subsection{2-term silting objects}

The study of support $\tau$-tilting pairs turns out to be equivalent to that of another class of objects: $2$-term silting objects~\cite[Section~3]{AdachiIyamaReiten}. Let $K^b(\proj \Lambda)$ be the homotopy category of bounded complexes of projective $\Lambda$-modules. Let $\cpx(\Lambda)$ be the full subcategory of $K^b(\proj \Lambda)$ consisting of \textit{$2$-term objects}, that is, complexes
\begin{gather*}
 \complexP = \cdots \to P_{m+1} \to P_m \to P_{m-1} \to \cdots
\end{gather*}
such that $P_m$ is zero unless $m\in \{0,1\}$. We will write $P_1\to P_0$ to denote the complex
\begin{gather*}
\cdots \to 0 \to P_1 \to P_0 \to 0 \to \cdots
\end{gather*}
A $2$-term object $\complexP$ is \textit{rigid} if $\Hom{K^b}(\complexP, \complexP[1]) = 0$. It is \textit{silting} if
\begin{itemize}\itemsep=0pt
 \item it is rigid, and
 \item $|\complexP| = |\Lambda|$.
\end{itemize}
This is a special case of a more general definition of silting objects, see~\cite{KellerVossieck}. Examples of $2$-term silting objects include $0\to \Lambda$ and $\Lambda \to 0$.

\begin{Definition} The \textit{$2$-term silting complex} of~$\Lambda$ is the simplicial complex~$\siltingComplex(\Lambda)$ whose vertices are isomorphism classes of indecomposable rigid $2$-term objects in $K^b(\proj \Lambda)$ and whose faces are sets of such objects whose direct sum is rigid. The \textit{exchange graph}~$\silt(\Lambda)$ is the dual graph of~$\siltingComplex(\Lambda)$, i.e., the graph whose vertices are isomorphism classes of basic $2$-term silting objects in $K^b(\proj \Lambda)$, and where two vertices are joined by an edge whenever the corresponding objects differ by exactly one direct summand.
\end{Definition}

For any $\Lambda$-module $M$, denote by $P_1^M\to P_0^M$ a minimal projective presentation of $M$.

\begin{Theorem}[{\cite[Theorem~3.2]{AdachiIyamaReiten}}] The map $(M,P) \mapsto \big(P_1^M\to P_0^M\big) \oplus (P\to 0)$ induces an isomor\-phism of simplicial complexes~$\tiltingComplex(\Lambda) \cong \siltingComplex(\Lambda)$, and thus of exchange graphs $\stau(\Lambda)$ $\cong \silt(\Lambda)$.
\end{Theorem}

\subsection{The $\b{g}$-vector of a 2-term object}

The results of this note rely on the following definition.

\begin{Definition} Let $\complexP$ be a $2$-term object in $\cpx(\Lambda)$. The \textit{$\b{g}$-vector} of $\complexP$, denoted by $\b{g}(\complexP)$, is the class of $\complexP$ in the Grothendieck group $K_0\big(K^b(\proj \Lambda)\big)$.
\end{Definition}
We will usually denote $\b{g}$-vectors as integer vectors by using the basis of the abelian group $K_0\big(K^b(\proj \Lambda)\big)$ given by the classes of the indecomposable projective modules $e_1 \Lambda, \dots, e_n \Lambda$ concentrated in degree $0$. Thus, if $\complexP$ is the $2$-term object
\begin{gather*}
 \bigoplus_{i \in [n]} (e_i \Lambda)^{\oplus b_i} \xrightarrow{} \bigoplus_{i \in [n]} (e_i \Lambda)^{\oplus a_i},
\end{gather*}
then its $\b{g}$-vector is $\b{g}(\complexP) = (a_i - b_i)_{i \in [n]}$.

In contrast to arbitrary $2$-term objects, rigid $2$-term objects are determined by their $\b{g}$-vector in the following sense.

\begin{Theorem}[{\cite[Sections~2.3 and 2.4]{DehyKeller}}]Let $\complexP$ and $\complexQ$ be two rigid $2$-term objects.
\begin{enumerate}[$(i)$]\itemsep=0pt
 \item If ${\b{g}(\complexP) = \b{g}(\complexQ)}$, then $\complexP$ and $\complexQ$ are isomorphic.
 \item The object $\complexP$ is isomorphic to an object of the form $(P_1\to P_0) \oplus \big(Q\stackrel{{\rm id}_Q}{\to} Q\big)$, where~$P_1$ and~$P_0$ do not have non-zero direct summands in common.
 \end{enumerate}
\end{Theorem}
Note that $\big(Q\stackrel{{\rm id}_Q}{\to} Q\big)$ is isomorphic to zero in $K^b(\proj \Lambda)$.

\section{Algebraic result}\label{sec:algRes}

We use the same notations as in the previous section. In particular, $\Lambda$ is a basic finite-dimensional $k$-algebra with complete set of pairwise orthogonal idempotents~$\{e_1, \dots, e_n\}$.

Let $J$ be a subset of $[n]$. We will study $2$-term objects that only involve the indecomposable projective modules $e_j \Lambda$ with $j\in J$.

\begin{Definition}Let $\cpx_J(\Lambda)$ be the full subcategory of $\cpx(\Lambda)$ whose objects are the $2$-term objects $P_1\to P_0$ such that all the indecomposable direct summands of $P_1$ and $P_0$ have the form~$e_j \Lambda$ with $j\in J$.
\end{Definition}

Our main interest will lie in the rigid objects in~$\cpx_J(\Lambda)$.

\begin{Definition}Let~$\siltingComplex_J(\Lambda)$ be the subcomplex of~$\siltingComplex(\Lambda)$ \textit{induced} by~$J$, that is, the subcomplex whose vertices are rigid objects in $\cpx_J(\Lambda)$.
Let $\silt_J(\Lambda)$ be the dual graph of~$\siltingComplex_J(\Lambda)$. Its vertices are isomorphism classes of basic objects $\complexP$ in $\cpx_J(\Lambda)$ satisfying
\begin{itemize}\itemsep=0pt
\item $\complexP$ is rigid;
\item if $\complexP'\in \cpx_J(\Lambda)$ and $\complexP\oplus \complexP'$ is rigid, then $\complexP'$ is a direct sum of direct summands of~$\complexP$.
\end{itemize}
Two vertices are joined by an edge whenever the corresponding objects differ by exactly one indecomposable direct summand.
\end{Definition}

In other words, the faces of $\siltingComplex_J(\Lambda)$ correspond to basic rigid objects whose $\b{g}$-vectors have zero coefficients in entries corresponding to elements not in $J$. In this sense, $\siltingComplex_J(\Lambda)$ is a~representa\-tion-theoretic analogue of the accordion complex~\cite{Chapoton-quadrangulations, GarverMcConville, MannevillePilaud-accordion} (see Theorem~\ref{thm:contractDiagonals}). This is the main motivation for the introduction of this object.

Let $e_J := \sum_{j\in J} e_j$, and consider the algebra $e_J \Lambda e_J$. Observe that $e_J \Lambda e_J$ is isomorphic to $\End{\Lambda}(\Lambda e_J)$. This has the following consequence. Let $\proj_J(\Lambda)$ be the full subcategory of $\proj(\Lambda)$ whose objects are isomorphic to direct sums of the indecomposable modules $e_j \Lambda$, with $j\in J$.

\begin{Lemma} The $k$-linear categories $\proj_J(\Lambda)$ and $\proj(e_J\Lambda e_J)$ are equivalent. In particular, the categories $K^b(\proj_J(\Lambda))$ and $K^b(\proj(e_J \Lambda e_J))$ are equivalent.
\end{Lemma}

This lemma immediately implies the following statement.

\begin{Theorem}The simplicial complexes~$\siltingComplex_J(\Lambda)$ and~$\siltingComplex(e_J \Lambda e_J)$ are isomorphic. In particular, their dual graphs $\silt_J(\Lambda)$ and $\silt(e_J \Lambda e_J)$ are isomorphic.
\end{Theorem}

\begin{Corollary}The simplicial complex~$\siltingComplex_J(\Lambda)$ is a pseudomanifold of dimension~$|J|-1$. In particular, its dual graph~$\silt_J(\Lambda)$ is $|J|$-regular.
\end{Corollary}

\section{Application: accordion complexes of dissections}\label{sec:accordions}

\begin{figure}[t]\centering \includegraphics[width=.9\textwidth]{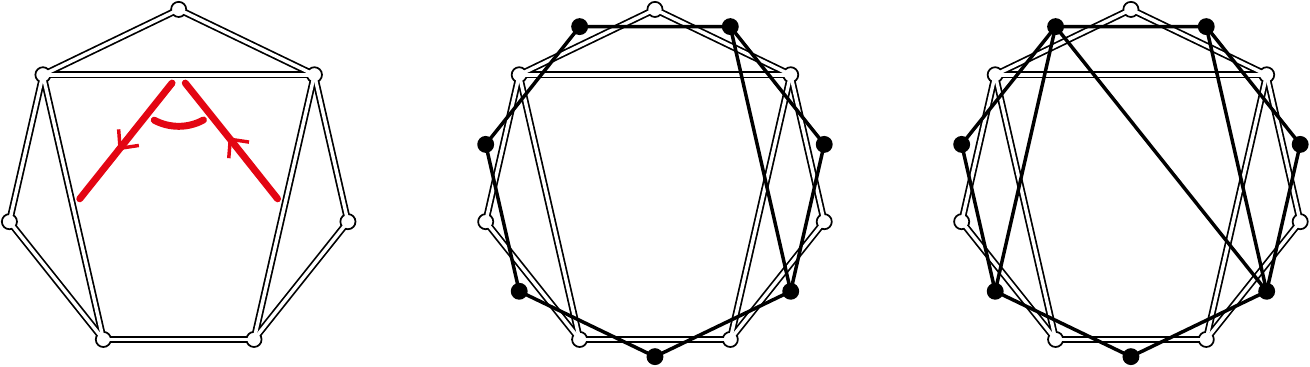}
	\caption{A dissection~$\dissection_\circ$ with its quiver~$\overline{\quiver}(\dissection_\circ)$ (left), a $\dissection_\circ$-accordion diagonal (middle) and a $\dissection_\circ$-accordion dissection (right).}
	\label{fig:exmAccordionDissections}
\end{figure}

Let~$\polygon$ be a convex polygon. We call \textit{diagonals} of~$\polygon$ the segments connecting two non-consecutive vertices of~$\polygon$. A \textit{dissection} of~$\polygon$ is a set~$\dissection$ of non-crossing diagonals. It dissects the polygon into \textit{cells}. We denote by~$\overline{\quiver}(\dissection)$ the quiver with relations whose vertices are the diagonals of~$\dissection$, whose arrows connect any two counterclockwise consecutive edges of a cell of~$\dissection$, and whose relations are given by triples of counterclockwise consecutive edges of a cell of~$\dissection$.
See Fig.~\ref{fig:exmAccordionDissections} for an example.

We now consider $2m$ points on the unit circle alternately colored black and white, and let~$\polygon_\circ$ (resp.~$\polygon_\bullet$) denote the convex hull of the white (resp.~black) points.
We fix an arbitrary reference dissection~$\dissection_\circ$ of~$\polygon_\circ$. A diagonal~$\delta_\bullet$ of~$\polygon_\bullet$ is a \textit{$\dissection_\circ$-accordion diagonal} if it crosses either none or two consecutive edges of any cell of~$\dissection_\circ$. In other words, the diagonals of~$\dissection_\circ$ crossed by~$\delta_\bullet$ together with the two boundary edges of~$\polygon_\circ$ crossed by~$\delta_\bullet$ form an accordion. A \textit{$\dissection_\circ$-accordion dissection} is a set of non-crossing \mbox{$\dissection_\circ$-accordion} diagonals. See Fig.~\ref{fig:exmAccordionDissections} for an example. We call \textit{$\dissection_\circ$-accordion complex} the simplicial complex~$\accordionComplex(\dissection_\circ)$ of $\dissection_\circ$-accordion dissections. This complex was studied in recent works of F.~Chapoton~\cite{Chapoton-quadrangulations}, A.~Garver and T.~McConville~\cite{GarverMcConville}, and T.~Manneville and V.~Pilaud~\cite{MannevillePilaud-accordion}.

For a diagonal~$\delta_\circ$ of~$\dissection_\circ$ and a $\dissection_\circ$-accordion diagonal~$\delta_\bullet$ intersecting~$\delta_\circ$, we consider the three edges (including~$\delta_\circ$) crossed by~$\delta_\bullet$ in the two cells of~$\dissection_\circ$ containing~$\delta_\circ$. We define~$\sign{\delta_\circ}{\dissection_\circ}{\delta_\bullet}$ to be $1$, $-1$, or~$0$ depending on whether these three edges form a~$\ZZZ$, a~$\SSS$, or a~$\VVV$. The \textit{$\b{g}$-vector} of~$\delta_\bullet$ with respect to~$\dissection_\circ$ is the vector~$\gvector{\dissection_\circ}{\delta_\bullet} \in \R^{\dissection_\circ}$ whose $\delta_\circ$-coordinate is~$\sign{\delta_\circ}{\dissection_\circ}{\delta_\bullet}$.

\begin{Example}\label{exm:associahedron}When the reference dissection~$\dissection_\circ$ is a triangulation of~$\polygon_\circ$, any diagonal of~$\polygon_\bullet$ is a $\dissection_\circ$-accordion diagonal. The $\dissection_\circ$-accordion complex is thus an $n$-dimensional associahedron (of type~$A$), where~$n = m-3$. As explained in \cite{CalderoChapotonSchiffler}, the $\dissection_\circ$-accordion complex is isomorphic to the $2$-term silting complex of the quiver~$\overline{\quiver}(\dissection_\circ)$ of the triangulation~$\dissection_\circ$. The~isomorphism sends a~diagonal of~$\polygon_\bullet$ to the $2$-term silting object with the same $\b{g}$-vector. See Fig.~\ref{fig:associahedron} for an illustration.
\end{Example}

With the notations we introduced, we can now restate Theorem~\ref{thm:mainDissection} more precisely.

\begin{Theorem}\label{thm:bijectionAccordionComplexSiltingComplex}For any reference dissection~$\dissection_\circ$, the $\dissection_\circ$-accordion complex is isomorphic to the $2$-term silting complex of the quiver~$\overline{\quiver}(\dissection_\circ)$.
\end{Theorem}

One possible approach to Theorem~\ref{thm:bijectionAccordionComplexSiltingComplex} would be to provide an explicit bijective map between $\dissection_\circ$-accordion diagonals and $2$-term silting objects for~$\overline{\quiver}(\dissection_\circ)$. Such a map is easy to guess using \mbox{$\b{g}$-vectors}, but the proof that it is actually a bijection and that it preserves compatibility is intricate. This approach was developed in the more general context of non-kissing complexes of gentle quivers with relations in~\cite[Propostion~2.44]{PaluPilaudPlamondon}. In this note, we use an alternative simpler strategy to obtain Theorem~\ref{thm:bijectionAccordionComplexSiltingComplex} by using Theorem~\ref{thm:mainAlgThm} and understanding accordion complexes as certain subcomplexes of the associahedron.

For that, consider two nested dissections~$\dissection_\circ \subset \dissection_\circ'$. Observe that any~$\dissection_\circ$-accordion diagonal is a $\dissection_\circ'$-accordion diagonal. Conversely a $\dissection_\circ'$-accordion diagonal~$\delta_\bullet$ is a $\dissection_\circ$-accordion diagonal if and only if it does not cross any diagonal~$\delta_\circ'$ of~$\dissection_\circ' \ssm \dissection_\circ$ as a $\ZZZ$ or a~$\SSS$, that is if and only if the~$\delta_\circ'$-coordinate of its $\b{g}$-vector~$\gvector{\dissection_\circ'}{\delta_\bullet}$ vanishes for any~$\delta_\circ' \in \dissection_\circ' \ssm \dissection_\circ$.
This observation shows the following statement.

\begin{Theorem}[{\cite[Section~4.2]{MannevillePilaud-accordion}}]\label{thm:contractDiagonals}For any two nested dissections~$\dissection_\circ \subset \dissection_\circ'$, the accordion complex~$\accordionComplex(\dissection_\circ)$ is isomorphic to the subcomplex of~$\accordionComplex(\dissection_\circ')$ induced by $\dissection_\circ'$-accordion diagonals~$\delta_\bullet$ whose $\b{g}$-vectors $\gvector{\dissection_\circ'}{\delta_\bullet}$ lie in the coordinate subspace spanned by elements in~$\dissection_\circ$.
\end{Theorem}

In order to prove Theorem~\ref{thm:bijectionAccordionComplexSiltingComplex} we now turn to associative algebras. Let~$\overline{\quiver} = (\quiver, \relations)$ be any gentle quiver with relations~\cite{ButlerRingel} and~$J$ be any subset of vertices of~$\overline{\quiver}$. We call \textit{shortcut quiver} the quiver with relations~$\overline{\quiver}_J = (\quiver_J, \relations_J)$ whose vertices are the elements of~$J$, whose arrows are the paths in~$\overline{\quiver}$ with endpoints in~$J$ but no internal vertex in~$J$, and whose relations are inherited from those of~$\overline{\quiver}$. Then the quotient~$k\quiver_J/\relations_J$ of the path algebra~$k\quiver_J$ is gentle and is isomorphic to the algebra~$e_J (k\quiver/\relations) e_J$.

\begin{Example}Quivers of dissections are shortcut quivers: if~$\dissection_\circ \subset \dissection_\circ'$, then~${\overline{\quiver}(\dissection_\circ) = \overline{\quiver}(\dissection_\circ')_{\dissection_\circ}}$.
In particular, for any dissection~$\dissection_\circ$, the quiver~$\overline{\quiver}(\dissection_\circ)$ is a shortcut quiver of the quiver with relations of a cluster tilted algebra of type~$A$.
\end{Example}

The following statement is an application of Theorem~\ref{thm:mainAlgThm} to gentle algebras.

\begin{Theorem}\label{thm:contractVertices}For any gentle quiver with relations~$\overline{\quiver}$ and any subset~$J$ of vertices of~$\overline{\quiver}$, the \mbox{$2$-term} silting complex~$\siltingComplex(\overline{\quiver}_J)$ for the shortcut quiver~$\overline{\quiver}_J$ is isomorphic to the subcomplex of the $2$-term silting complex~$\siltingComplex(\overline{\quiver})$ induced by $2$-term silting objects whose $\b{g}$-vectors lie in the coordinate subspace spanned by vertices in~$J$.
\end{Theorem}

Combining Theorems~\ref{thm:contractDiagonals} and~\ref{thm:contractVertices} together with Example~\ref{exm:associahedron}, we obtain Theorem~\ref{thm:bijectionAccordionComplexSiltingComplex} (and Theorem~\ref{thm:mainDissection}).

\section{Concluding remarks}

\begin{Remark}There is a geometric interpretation of the common phenomenon described in Theorems~\ref{thm:contractDiagonals} and~\ref{thm:contractVertices}. For a $\dissection_\circ$-accordion dissection~$\dissection_\bullet$, denote by~$\R_{\ge0}\,\gvector{\dissection_\circ}{\dissection_\bullet}$ the polyhedral cone generated by the set of $\b{g}$-vectors~$\gvector{\dissection_\circ}{\dissection_\bullet} := \set{\gvector{\dissection_\circ}{\delta_\bullet}}{\delta_\bullet \in \dissection_\bullet}$. The collection~$\gvectorFan$ of cones~$\R_{\ge0}\,\gvector{\dissection_\circ}{\dissection_\bullet}$ for all $\dissection_\circ$-accordion dissections~$\dissection_\bullet$ is a complete simplicial fan called \mbox{\textit{$\b{g}$-vector fan}} of~$\dissection_\circ$~\cite{MannevillePilaud-accordion}. The crucial feature of this fan is that no coordinate hyperplane meets the interior of any of its maximal cones. This is often referred to as the \textit{sign-coherence property} of $\b{g}$-vectors. It implies that for any two nested dissections~$\dissection_\circ \subset \dissection_\circ'$, the section of~$\gvectorFan[\dissection_\circ']$ with the coordinate subspace~$\R^{\dissection_\circ}$ is a subfan of~$\gvectorFan[\dissection_\circ']$. The content of Theorem~\ref{thm:contractDiagonals} is that this subfan is the $\b{g}$-vector fan~$\gvectorFan$. A~similar statement holds for Theorem~\ref{thm:contractVertices}.
\end{Remark}

\begin{Remark}In the theory of cluster algebras, a standard operation consists of freezing a subset of the initial cluster. This corresponds to taking a section of the $\b{d}$-vector fan by a coordinate subspace. To the best of our knowledge, the same operation on the $\b{g}$-vector fan studied in this note was not considered before in the literature.
\end{Remark}

\begin{Remark} The connection between representation theory and accordion complexes was already considered by A.~Garver and T.~McConville in~\cite[Section~8]{GarverMcConville}. However, their approach deals with $\b{c}$-vectors and simple-minded collections while our approach deals with $\b{g}$-vectors and silting objects.
\end{Remark}

\subsection*{Acknowledgements}

We are grateful to F.~Chapoton for pointing out to us the isomorphism between the two graphs of Fig.~\ref{fig:accordiohedron} which gave us the motivation for the present note. We also thank R.~Schiffler for his comments on a previous version. Finally, we are grateful to an anonymous referee for helpful suggestions on the presentation of this note. The first two authors are partially supported by the French ANR grant SC3A~(15\,CE40\,0004\,01). The last author is supported by the ISF grant 1144/16.

\pdfbookmark[1]{References}{ref}
\LastPageEnding

\end{document}